\documentclass{amsart}
\usepackage{times,hyperref,inputenc,graphicx}
\newtheorem{theorem}{Theorem}[section]

\theoremstyle{definition}
\newtheorem{definition}[theorem]{Definition}
\theoremstyle{remark}
\newtheorem{remark}[theorem]{Remark}
\numberwithin{equation}{section}
\newcommand{\be}{\begin{equation}}
\newcommand{\ee}{\end{equation}}
\newcommand{\bd}{\begin{displaymath}}
\newcommand{\ed}{\end{displaymath}}
\newcommand{\bea}{\begin{eqnarray}}
\newcommand{\eea}{\end{eqnarray}}
\newcommand{\beas}{\begin{eqnarray*}}
\newcommand{\eeas}{\end{eqnarray*}}
\newcommand{\bc}{\begin{center}}
\newcommand{\ec}{\end{center}}

\title[Warped-like product  manifolds with exceptional holonomy groups]
{Warped-like product  manifolds with exceptional holonomy groups}%
\author{Selman Uguz}
\address{ Selman Uguz, UT DALLAS, DEPT. MATH. SCIENCES, RICHARDSON, TX 75080} \email{{\tt  selman.uguz@utdallas.edu}}

\date{}%
\thanks{}
\begin{document}
\maketitle
\begin{abstract}
In this paper we review   $G_2$ and $Spin(7)$ geometries in
relation with
  a special type of   metric structure
which we call  warped-like product metric. We present a general
ansatz of warped-like product metric as a definition of
warped-like product. Considering fiber-base decomposition, the
definition of warped-like product is regarded as a generalization
of multiply-warped product manifolds, by allowing the fiber metric
to be non block diagonal. For some special cases, we present
explicit example of $(3+3+2)$  warped-like product manifolds with
$Spin(7)$ holonomy of the form $M=F\times B$, where the base $B$
is a  two dimensional Riemannian manifold, and the fibre $F$  is
of the form $F=F_1\times F_2$ where $F_i$'s $(i=1,2)$ are
Riemannian $3$-manifolds. Additionally an explicit example of
$(3+3+1)$ warped-like product manifold with $G_2$ holonomy  is
studied. From the literature, some other special warped-like
product metrics with $G_2$ holonomy are also presented in the
present study. We believe that our approach of the warped-like
product metrics will be an important notion for the  geometries
which use warped and multiply-warped product structures, and
especially manifolds with exceptional holonomy.


\end{abstract}

\section{Introduction}

A warped product metric   is an important notion in Riemannian
geometry as well as in physics and related areas. %How it curves space will define
%one or other solution to a space-time world.
Using its properties, the different works of geometry and physics
use warped product geometries. Furthermore many basic solutions of
the Einstein field equations are used of warped geometries, e.g.
the Schwarzschild solution and the Robertson-Walker models
\cite{O'Neil83}. In string theory, Yau  in \cite{Yau}
  discussed that ``{
...the easiest way to partition the ten-dimensional space is to
cut it cleanly, splitting it into four-dimensional space-time and
six-dimensional hidden subspace... and in the non-kahler case, the
ten-dimensional spacetime is not a Cartesian product but rather} {
\it a warped product}." However, the fundamental rigidity theorems
for manifolds of non-negative or positive Ricci curvature are the
volume cone  (metric cone) theorem, the maximal diameter theorem
and the splitting theorem \cite{Cheg}. Each theorem asserts that
if a certain geometric value such as volume or diameter  is
possible large enough relative to the pertinent lower bound on
Ricci curvature, then the metric of the manifold  is a {\it warped
product metric} of a special type (see details in \cite{Cheg}).

On the other hand, the notion of holonomy group of a Riemannian
manifold defined by \'{E}lie Cartan   has exhibited to be an
efficient tool in the study of Riemannian manifolds  (see for the
details \cite {KN,Ilka,SW}).
 The list of all possible holonomy groups of irreducible, simply-connected,
non-symmetric spaces was given by M. Berger in 1955 \cite{Berger}.
Berger's list (refined  by Alekseevskii \cite{A1} and Gray-Brown
\cite{GrayB}) includes the groups $SO(n)$ in $n$-dimensions,
$U(n)$, $SU(n)$ in $2n$-dimensions, \ $Sp(n)$,
 $ Sp(n)Sp(1)$ in $4n$-dimensions  and two special cases, $G_2$
holonomy in $7$-dimensions and $Spin(7)$ holonomy  in
$8$-dimensions. Manifolds with holonomy $SO(n)$ constitute the
generic case, all others are denoted as manifolds with ``{\it
special holonomy}" and the last two cases are described as
manifolds with ``{\it exceptional holonomy}".

 The existence  of manifolds with exceptional holonomy was first demonstrated
 by R.Bryant \cite{R}, complete examples were given by  R. Bryant and S. Salamon {\cite{RS}} and
the first compact examples were found by  D. Joyce in 1996
\cite{J}. The study of manifolds with exceptional holonomy  and
the construction of explicit examples with special types (e.g.
Fernandez classifications \cite{Fer})
 are still an active research area in mathematics
and physics (see also references in \cite {SW,U,U3,
Fer,Sema,Cab}).

In physics, there exist a special interest in the construction of
$G_2$ and $Spin(7)$ holonomy metrics
 due to their
application in supergravity compactification %for preserving
%certain amount of supersymmetry
(see more details \cite{Ilka,R,G1,Jb,Sema}). Since manifolds with
special holonomy provide some geometrical structures for reducing
the number of super-symmetries, they are natural candidates for
the extra dimensions in string and M-theory \cite{Ilka,F11,F12,
GPP,G1, Jb, Sema}. Following the constructions of Bryant and
Salamon \cite{R,RS} and   Joyce \cite{Jb}, the study of gauge
fields on reduced holonomy manifolds was formulated in the physics
and mathematics literature. Existence theorems have been given on
both compact and non-compact spaces (see refs. in \cite{Clarke}).

In \cite{Baza}, complete Riemannian metrics with holonomy group
$G_2$ on manifolds obtained by deformation of cones over $S^3
\times S^3$ are constructed in details. Their   idea of the paper
was also used for the purpose of constructing complete Riemannian
metrics with holonomy group $Spin(7)$ (see details in
\cite{Baza}). %It is summarized as follows, consider the standard
%conic metric over a Riemannian manifold with special geometry.
%Then any deformation of this metric depends on a number of
%functional parameters using which the $G_2 \ (or \ Spin(7))$
%structure can be explicitly specified in this study (for the
%details in \cite{Baza}).
In \cite{Clarke}, it is given a construction of $G_2$ and
$Spin(7)$ instantons on exceptional holonomy manifolds constructed
by Bryant and Salamon \cite{RS}, by using an ansatz of spherical
symmetry coming from the manifolds being the total spaces of
rank-4 vector bundles. For the seven dimensional $G_2$ case, it is
shown that the connections are asymptotic to Hermitian Yang-Mills
connections on the nearly Kahler $S^3 \times S^3$   in the
asymptotically conical model \cite{Clarke}.

The motivation for this work was the results \cite{BU,U,U2} and
the explicit $Spin(7)$ metric special solution on $M=F_1\times F_2
\times B=S^3\times S^3\times R^2$  in \cite{YO}, with the
uniqueness of this solution proved in \cite{BU}.  It is noticed
that their metric ansatz (structure) was a generalization of
warped products and we presented   ``(3+3+2) warped-like product
metrics" as a general framework for the special metrical ansatz as
given in \cite{BU}. It is looked whether one could obtain other
solutions by relaxing some of their assumptions, in particular
without requiring the three dimensional submanifolds to be $S^3$.
It is also proved that the connection of the fibers is determined
by the Bonan form $\Omega $ given \cite{BU}. Considering suitable
global assumptions, it is concluded that the fibers  are
$3$-spheres with constant curvature $k>0$. Similar computation is
done for 7-dimensional manifold with $G_2$ holonomy as a $(3+3+1)$
special warped-like product manifold \cite{U}.

In the present paper, we   study   more general   metric structure
called ``warped-like product metric". It is presented a general
ansatz of warped-like product metric as a definition of
warped-like product (see Definition 4.2). Using fiber-base
decomposition, the definition of warped-like product is considered
as a generalization of multiply-warped product manifolds, by
allowing the fiber metric to be non block diagonal. For some
special cases, we present explicit examples warped-like product
manifolds with $G_2$ and $Spin(7)$ holonomies (i.e. exceptional
holonomy).

The outline of paper is given as follows: we set up the basic
classifications of $G_2$ and $Spin(7)$ geometries in Section 2.
Warped product and a generalization of warped products which is
called multiply-warped product \cite{FS} are presented in Section
3.
  In Section 4,  we  study a
generalization of multiply-warped product manifolds as a
``warped-like product", by allowing the fiber metric to be non
block diagonal.    Some of the explicit examples of warped-like
product manifolds with $G_2$ and $Spin(7)$ holonomies are studied
in Section 5.
 The conclusions and planning of the next studies are presented in Section 6.

\section{Exceptional geometries in special dimensions $d=7$
and $d=8$} In this section we present the basics of $G_2$ and
$Spin(7)$ geometries with their related properties
\cite{Ilka,R,Jb,Salamonb,Fer,Iva1,Iva2}. Let us start with seven
dimensional case. Suppose that $N$ indicates a $7$-dimensional
manifold with an invariant $G_2$ structure. Thus, $N$ is endowed
with  a non-degenerate $3$-form $\varphi$ that induce a Riemannian
metric. The fundamental material for the $G_2$ (also $Spin(7)$)
geometry can be found in standard holonomy references books
~\cite{Jb,Salamonb}. Let us only recall that the Riemannian
geometry of $N$ is completely determined by the special form
(called Fundamental $3$-form in $7$-dimension)
\begin{equation} \label{34}
     \varphi=e^{125}-e^{345}+e^{567}+e^{136}+e^{246}-e^{237}+e^{147}.
\end{equation}
It has become customary to suppress wedge signs when writing
differential forms, so $e^{ij\dots}$
  indicates $e^i\wedge e^j\wedge\dots$ from now on.
The results of Fernandez and Gray~\cite{Fer1} give one to describe
$G_2$ geometry exclusively in algebraic terms, by looking at the
various components of $d\varphi, d*\varphi$ in the irreducible
summands. Many authors have studied special classes of $G_2$
structures, see for instance~\cite{F11,F12,Sema,Iva1,Cab1}. Note
that 16 classes of $G_2$ manifolds in the Fernandez
    classification can be described in terms of the Lee form
    summarized some of them as given in Table 1.

{\small  $$ \begin{tabular}{|c|c|}
  \hline
 \qquad  {\bf Some classes  of $G_2$ manifolds} \qquad  & \qquad   {\bf Conditions} \qquad  \\
 \hline
   Parallel case:  & $d\varphi=0$ and $d*\varphi=0$  \\  \hline
   Almost Parallel or closed  (or calibrated symplectic ) case:  & $d\varphi=0$
   \\
    \hline
 Locally conformally (almost) parallel case:  & $d\varphi= \frac{3}{4}\Theta \wedge
   \varphi$ \\
 \hline
   Balanced case  & $\delta\varphi=0$ and $d \varphi\wedge \varphi=0$  \\  \hline
 Locally conformally (almost) parallel case:  & $d\varphi= \frac{3}{4}\Theta \wedge
   \varphi$ \\ \hline
   Locally conformally parallel case:  &$d\varphi= \frac{3}{4}\Theta \wedge
   \varphi$ and $d*\varphi= \Theta \wedge
   *\varphi$
   \\ \hline
  Cocalibrated   (or semi-paralel or cosymplectic) case& $\delta\varphi=0$   \\
  \hline
\end{tabular} $$
  \begin{center} Table 1:   Special types of some $G_2$
structures in dimension seven.
\end{center} \vspace{2mm}
}

  Moving up one dimension,
we consider a product $M$ of $N$ with $\mathbb{R}$, endowed with
metric $g$. Indicating by $e^8$ the unit 1-form on the real line
one obtains a basis for the cotangent spaces $T_p^*M$.  The
manifold $M$ inherits a non-degenerate four-form
\begin{equation} \label{ww}\Omega=\varphi\wedge e^8+*\varphi
\end{equation} which %is stable differential form (in the sense of
%Hitchin \cite{Hitchin:forms}, it occurs in some special
%dimensions)
%    and
defines a reduction to the Lie group \cite{Iva1}. In equation
(\ref{ww}), a special note that the Hodge dual map of $\varphi$
(i.e. $*\varphi$) is considered on $7$-dimensional manifold $N$,
as in the equation (\ref{34}) \cite{Jb}.

In $Spin(7)$ geometry,  we recall that the special $4$-form
$\Omega$ is self-dual $*\Omega=\Omega$, where $*$ is the Hodge
operator and the $8$-form $*\Omega\wedge \Omega$ coincides with
the volume form. It is well known that the subgroup of $GL(8,R)$
which fixes $\Omega$ is isomorphic to the double covering
$Spin(7)$ of $SO(7)$ \cite{Salamonb}. Moreover, $Spin(7)$ is a
compact simply-connected Lie group of dimension $21$ \cite{Ilka}.

The $4$-form $\Omega$ corresponds to a real spinor $\phi$ and
therefore, $Spin(7)$ can be identified as the isotropy group of a
non-trivial real spinor \cite{Cab}. A $3$-fold vector cross
product $P$ on $R^8$ can be defined by \begin{equation}< P(x\wedge
y u\wedge z), t
>= \Omega(x, y, z, t), \ for \ \ x, y, z, t \in R^8.\end{equation} Then
$Spin(7)$ is also characterized by \begin{equation}Spin(7) = \{a
u\in O(8)|P(ax u\wedge ay u\wedge az) = P(x u\wedge y u\wedge z),
x, y, z u\in R^8 \}.\end{equation} The inner product $<,>$ on
$R^8$ can be reconstructed from $\Omega$ \cite{Fer, Cab1}, which
corresponds with the fact that $Spin(7)$ is a subgroup of $SO(8)$.
A $Spin(7)$ structure on an $8$-dimensional manifold $M$ is by
definition a reduction of the structure group of the tangent
bundle to $Spin(7)$, we shall also say that $M$ is a $Spin(7)$
manifold. This can be described geometrically by saying that there
is a $3$-fold vector cross product $P$ \cite{Fer} defined on $M$,
or equivalently there exists a nowhere vanishing differential
$4$-form $\Omega$ on $M$ which can be locally written as
\begin{eqnarray}
\Omega &=& e^{1258} + e^{3458} + e^{1368} - e^{2468} + e^{1478} +
e^{2378} - e^{5678} \nonumber \\
{}&{}& - e^{1267} - e^{3467}+ e^{1357} - e^{2457} - e^{1456} -
e^{2356} + e^{1234}.
\end{eqnarray}
This special $4$-form $\Omega$ is called ``Bonan (Cayley or
Fundamental) form" of the $Spin(7)$ manifold $M$
\cite{R,Bonan,BU}. We also recall that a $Spin(7)$ manifold $(M,g,
\Omega )$ is said to be parallel if the holonomy of the metric
$Hol(g)$ is a subgroup of $Spin(7)$.
 This is equivalent to
saying that the fundamental form $\Omega$ is parallel with respect
to the Levi-Civita connection u$\nabla^{LC}$ of the metric g.
Moreover, $Hol(g) u\subset Spin(7)$ if and only if $d\Omega = 0$
\cite{R} and any parallel $Spin(7)$ manifold is Ricci-flat
\cite{Bonan}.

According to the classification given by   Fernandez \cite{Fer},
there are four classes of $Spin(7)$ manifolds obtained as
irreducible representations of $Spin(7)$ of the space
u$\nabla^{LC}\Omega$. By using the fact that given by Cabrera et.
al. \cite{Cab1}, It is considered the $1$-form   of the
$8$-manifold defined by
\begin{equation} 7\Theta =-*(*d\Omega \wedge \Omega ) = *(\delta
\Omega \wedge \Omega).\end{equation} It is called the Lee form
(this $1$-form is denoted by $\Theta$) of a given $Spin(7)$
structure \cite{Iva}. The four classes of $Spin(7)$ manifolds in
the Fernandez classification can be described in terms of the Lee
form as below \cite{U,Iva}: $\label{ww}W_0 : d\Omega = 0; \ \ W_1
: \Theta = 0; \ \ W_2 : d\Omega = \Theta \wedge \Omega  ; \ \ W_4
= W_1 \oplus W_2. $
%A Spin(7) structure of class $W_1$, that is, a $Spin(7)$ structure
%with Lee form equal to zero, is called a balanced $Spin(7)$
%structure
%  Cabrera \cite{Cab} shows that the Lee form of a
%$Spin(7)$ structure in the class $W_2$ is closed and therefore
%such a manifold is locally conformally equivalent to a parallel
% $Spin(7)$  manifold and it is called locally conformally parallel. If
%the Lee form is not exact (i.e. the structure is not globally
%conformally parallel), it is called strict locally conformally
%parallel.
We summarize these facts in the following Table 2.

$$\begin{tabular}{|c|c|}
  \hline
 \qquad  {\bf The classes  of Spin(7) manifolds} \qquad  & \qquad   {\bf Conditions} \qquad  \\
  \hline
   Parallel case: $W_0$ & $d\Omega=0,\Theta=0$ \\  \hline
  Balanced case: $W_1$&$\Theta=0$ \\ \hline
   Locally conformally parallel case: $W_2$&$d\Omega= \Theta \wedge \Omega$
   \\ \hline
   Mixed type: $W_4=W_1+W_2$& - \\
  \hline
\end{tabular} \\   $$
 \begin{center} Table 2:   Fernandez classification
table of $Spin(7)$ manifolds. \end{center} \vspace{2mm}

\section{Multiply-warped  and warped product manifolds}

We start with the definition of   multiply-warped and warped
product manifolds. In the next subsection, we present the
definition of our special chosen warped product; warped-like
product manifolds as a generalization of a multiply-warped
product.

\subsection{Warped  product manifolds}

Let $(F,g_F)$, $(B,g_B)$ be Riemannian manifolds and $f
>0$  be smooth function on $B$. A  {\it warped product manifold } is a
product manifold $M=F\times B$ equipped with the metric
 \be g= \pi_2^* g_B+ (f\circ \pi_2)^2\pi_1^* g_F, \ee
where   $\pi_1: F \times B \longrightarrow F$ and $\pi_2: F \times
B \longrightarrow B$ are the natural projections. \be g=%\left(%
%\begin{array}{c|c}
%  g_F & 0 \\ \hline
%  0 &g_B \\
%\end{array}%
%\right) =
\left(%
\begin{array}{c|c}
  fg_{F} & 0  \\ \hline 0  & g_B\\
\end{array}%
\right) \ee

That is, in local coordinates the first block that depends on the
coordinates of the first group of coordinates is multiplied by a
function of the second group of coordinates. If the definition
holds  an open subset of $M$, then $M$ is called {\it locally
warped product} manifold. Basic properties of warped product
manifolds can be found  in \cite{O'Neil83}.

\subsection{Multiply-warped product manifolds}

A generalization of the notion of warped product metrics  is the
``multiply-warped products", defined as follows \cite{FS}.
 Let $(F_i,g_i)$,
$i=1,2,...,k$ and   $(B,g_B)$ be Riemannian manifolds and  $f_i
>0$  be smooth functions on $B$. A {\it multiply-warped product manifold }
 is the product manifold \be F_1 \times
F_2 \times ...\times F_k \times B,\ee  equipped with the metric
\be \label{mul} g=\pi_B^* g_B +\sum_{i=1}^k(f_i\circ \pi_B)^2
\pi_i^* g_i, \ee where   $\pi_B: F_1 \times F_2 \times ...\times
F_k \times B \longrightarrow B$ and $\pi_i: F_1 \times F_2 \times
...\times F_k \times B \longrightarrow F_i$ are the natural
projections on $B$ and $F_i$ respectively.  In this scheme, the
metric is block diagonal, with the metrics of the $F_i$'s  are
multiplied by a conformal factor depending on the coordinates of
the base.

We shall first give a local  description of the generalization of
the warped product structure.  We shall use the summation
convention whenever appropriate. The base-fiber decomposition
suggests that the local formulation of the line element is of the
form \be ds^2=(g_F)_{ij}^{ab}\ dy_a^i\otimes dy_b^j +  (g_B)_{ij}
dx^i\otimes dx^j,$$ where $g_B$ depends on the base coordinates
and $g_F$ depends both on the base and the fiber coordinates. For
a  generalized warped product the local coordinate expression of
$g_F$ is of the form
$$g=\left(%
\begin{array}{c|c}
  g_F & 0 \\ \hline
  0 &g_B \\
\end{array}%
\right) =\left(%
\begin{array}{cccccc|c}
  A_1g_{1} & 0 & 0 & 0 & 0 & 0& 0 \\
  0 & A_2g_{2} & 0 & 0 & 0 & 0 & 0\\
  0 & 0 & A_3g_{3} & 0 & 0 & 0 & 0\\
  \dots & \dots & \dots & \dots & \dots & \dots & 0\\
  \dots & \dots & \dots & \dots & \dots & \dots & 0\\
  0 & 0 & 0 & 0 & 0 & A_kg_{k} & 0\\ \hline 0 & 0 & 0 & 0 & 0 &   & g_B\\
\end{array}%
\right) \ee where the $A_a$'s are functions of the base
coordinates and each $g_a$ the metric on $F_a$ hence is a
symmetric matrix whose entries are functions of the coordinates of
the fiber $F_a$.

 \vskip 0.2cm

We will study a special extension of multiply warped product type
warped metric construction  case called {\it
 warped-like }
product
manifold in the following.%  in the following section.

\section {Warped-like product  manifolds}

 Let $M$   be the topologically
 product manifold
\be M=F_1\times F_2 \times ...\times F_k \times B, \ee
 where $dim \ F_a=n_a, \ (a=1,...,k)$ and $dim \ B=n$.
The manifold $B$ will be called the ``base" and the manifolds
$F_a$ will be called ``fibers".
 We shall
denote the local coordinates on the base $B$ by
$$\{x^1,\dots,x^n\}$$
and the local coordinates each fiber $F_a$ by
$$\{y_a^1,\dots,y_a^{n_a}\}.
$$

We shall work with Riemannian metrics and use the the summation
convention whenever appropriate.
%If the local coordinates on an $n$-dimensional  manifold are $\{t_i\}$, the  coordinate expression of the metric is given by
%$$ds^2=g_{ij}dt_i\otimes dt_j.$$
A local ``moving frame" on an $n$-dimensional Riemannian manifold
consists of $n$ local orthonormal sections of its  cotangent
bundle of. If we denote these local sections by $e_i$,  the metric
is given by
$$ds^2=e_1\otimes e_1+ e_2\otimes e_2+\dots+ e_n\otimes e_n.$$
Thus, the metric written with respect to the moving frame has
constant coefficients.  The moving frame has redundancies; we may
transform the frame as
$$e'_i=P_{ij}e_j$$
where $P$ belongs to $SO(n)$; such transformations leave the
metric invariant and they are called the ``gauge transformations".

The ``warped-like" manifold structure that we are trying to define
aims to express the metric with respect to a frame in such a way
that the functions that appear in the metric depend only on the
coordinates of the base.

To start with, let us assume that each manifold $F_a$ is equipped
with a Riemannian metric. Hence we have local orthonormal frames
on each of them. We denote these local orthonormal frames  by
$\{\theta_a^i\}$ for $a=1,\dots ,k$ and $i=1,\dots, n_a$. Then
\begin{equation}\label{22}\{\theta_1^1,\dots,\theta_1^{n_1},\theta_2^1,\dots,\theta_2^{n_2},\dots,\theta_k^1,\dots,\theta_k^{n_k}
,dx^1,\dots,dx^n\} \end{equation} are linearly independent local
sections of  $T^*M$. Note the notation is chosen such that in
$\theta_a^i$, the lower index indicates the frame while the upper
index runs through the dimension of the fiber. We can define the
metric of the manifold $M$ by defining linearly independent local
sections of $T^*M$ and declaring these lo sections to be
orthonormal. We shall define a ``locally warped-like product"
metric on $M$ by declaring the following local sections to be
orthonormal: \be e_a^i=\sum_{b=1}^{k} \sum_{j=1}^{n_b}
A_{aj}^{ib}(x^l) \theta_b^j,\quad a=1,\dots, k, \quad i=1,\dots,
n_a.  \label{+++} \ee
$$e_\alpha=\sum_{\beta=1}^{n} A_\alpha^\beta dx^\beta.$$
Thus,  the metric written with respect to (\ref{22}) has
coefficients depending on the base coordinates only.

\begin{remark} The form of the $e_a^i$ given above is non
invariant under local orthonormal transformations of the fibers.
If
$$e_a^i\to \sum_{l=1}^{n_b} P_{al}(y_a^k) \theta_a^l,$$
then the moving frame transforms as \be e_a^i=\sum_{b=1}^{k}
\sum_{j=1}^{n_b} \sum_{l=1}^{n_b} A_{aj}^{ib}(x^l) P_{bl}(y_b^k)
\theta_b^l,\quad a=1,\dots, k, \quad i=1,\dots, n_a.\label{++} \ee
Therefore we should define a manifold to be warped-like if it
admits a moving frame of the form (\ref{+++}). \end{remark}

\vskip 0.2cm \noindent The functions that appear in (\ref{++})
that are sums of products of two different groups of coordinates
are called ``separable functions". %But these coefficients are not
%arbitrary separable functions, hence we cannot give a definition
%stating that the coefficients are separable functions.

\begin{definition} Let $M$   be the
topologically
 product manifold
\be M=F_1\times F_2 \times ...\times F_k \times B, \ee
 where $dim \ F_a=n_a, \ (a=1,...,k)$ and $dim \ B=n$.
Let the local coordinates on the base $B$ be $\{x^1,\dots,x^n\}$.
A metric on $M$ is called to be locally warped-like product, if
there are local frames $\{\theta_a^i\}$ on each $F_a$ that are
orthonormal with respect to the metric of the $F_a$, such that \be
e_a^i=\sum_{b=1}^{k} \sum_{j=1}^{n_b} A_{aj}^{ib}(x^l)
\theta_b^j,\quad a=1,\dots, k, \quad i=1,\dots, n_a. \label{+} \ee
$$e_\alpha=\sum_{\beta=1}^{n} A_\alpha^\beta dx^\beta$$
is an orthonormal frame for $M$. Note that \bea A_{aj}^{ib}=
A_{aj}^{ib}(x_1,x_2,...,x_n),
  A_{Bj}^{i}= A_{\alpha}^{\beta}(x_1,x_2,...,x_n). \eea \end{definition} \hfill $\square$

%\end{document}

% The local expression of the fiber space metric of generalized warped product  with respect to
%this basis is
%$$g_F=\left(%
%\begin{array}{ccccccc}
%  A_1I_{n_1} & 0 & 0 & 0 & \dots & 0 &0\\
%  0 & A_2I_{n_2} & 0 & 0 & \dots & 0 &0\\
%  0 & 0 & A_3I_{n_3} & 0 & \dots & 0 &0\\
%  \dots & \dots & \dots & \dots & \dots & \dots&0 \\
%  \dots & \dots & \dots & \dots & \dots & \dots&0 \\
%  0 & 0 & 0 & 0 & 0 & A_{k-1}I_{n_{k-1}}&0 \\
%   0 & 0 & 0 & 0 & 0 & 0&A_kI_{n_k} \\
%\end{array}%
%\right),
%$$
%where $I_{n_a}$ denotes the identity matrix of size $n_a$ and the
%$A_a$'s are functions of the base coordinates only.

We are looking for a generalization of the form where the local
coordinate expression of the metric is no longer block diagonal,
but still with a certain simple structure. Since we want to have a
fiber structure, we shall allow dependencies on the coordinates of
the base everywhere.  We shall discuss possible generalizations
below.

\begin{itemize}
    \item One possibility is to have each entry of $g_F$ to
        depend on all fiber coordinates, but be of the form
$$(g_F)_{ij}=\phi(x^1,\dots,x^n)\psi(y_1^1,\dots,y_k^{n_k}).$$
This means that there is a fiber-base decomposition, but no
internal decomposition of the fiber.
    \item Another possibility is to partition  $g_F$ according
        to the dimensions of the $F_a$'s. If we denote the
        block in this partition as $A_{ab}$ where $A_{ab}$ is
        a matrix with $n_a$ rows and $n_b$ columns, then we
        may require that
    $$(A_{ab})_{ij}=\phi_{ij}(x^1,\dots,x^n)\psi_{ij}(y_a^1,\dots,y_a^{n_a},y_b^1,\dots,y_b^{n_b}).$$
    \item A stronger form of this is to require that
    $$(A_{ab})_{ij}=\phi(x^1,\dots,x^n)(\psi_{ab})_{ij}(y_a^1,\dots,y_a^{n_a},y_b^1,\dots,y_b^{n_b}),$$
    that is to require that the same function of the base
    coordinates multiplies the whole submatrix.
 \end{itemize}

In the case of generalized warped products, the local expression
of the metric is as in case (3) with only diagonal blocks.  Then
each block is the metric of the fiber multiplied by a function of
the base coordinates and we can write the metric in terms of
pull-backs of the metrics of the components.

In all of the schemes above, above a fixed point $\{x^i\}$ of the
base we have the same situation.  That is for each $a$, $A_{aa}$
is a metric on the fiber component $F_a$ and similarly for each
pair $(a,b)$ the submatrix
\be g_{ab}=\left(%
\begin{array}{cc}
  A_{aa} & A_{ab} \\
  A_{ab} & A_{bb} \\
\end{array}%
\right) \ee is the metric on $F_a\times F_b$ (above this point).
But if we want to write the metric as a pull back, we should use
the third scheme and in addition  the off-diagonal block $A_{ab}$
should be a ``well defined" tensor on the product manifold
$F_a\times F_b$.
 Under these conditions we can write the metric on $M$ in terms of
 various pull-backs.  If $g_a$ be the metric on the $F_a$ and
 $h_{ab}$ the tensor corresponding to the off-diagonal block, then
\be g=\pi_B^*g_B+\sum_{a=1}^k \pi_a^*g_a+ \sum_{a<b=1}^k
\pi_{ab}^*h_{ab}.\ee Thus the main problem is the characterization
of $h_{ab}$ (see for special $3+3+2$ case in \eqref{3219} and
\eqref{mulpwarp33}). Our future research is on towards that
direction.

We shall now study the construction of these metrics in terms of
an orthonormal frame. If $\{e_a\}$, $a=1,\dots ,N$ is an
orthonormal frame  then the metric is
$$ds^2=e_1\otimes e_1+\dots +e_N\otimes e_N.$$
Denoting by ${\bf e}$ the vector whose entries are the 1-forms
$e_i$ and omitting the tensor product sign, we can write this as
$$ds^2 = {\bf e}^t {\bf e}.$$
If the frame transforms as
$${\bf e}= P\tilde{\bf e},$$
then
$$ds^2 = (P\tilde {\bf e})^t P\tilde{\bf e}.$$
Thus if we have the metric to be invariant we should have
$$P^t P=I,$$
that is the allowable frame rotations should belong to the
orthogonal group if an arbitrary Riemannian metric has to be
preserved.  This is the reduction of the structure group from
$Gl(N,R)$ to $O(N,R)$.

%
%
%In all cases, there is
%
%$$g_F=\left(%
%\begin{array}{cccccc}
%  A_1g_{1} & B_{12}h_{12} & B_{13}h_{13} & \dots & \dots & 0 \\
%  B_{21}h_{21} & A_2g_{n_2} & 0 & 0 & 0 & 0 \\
%  0 & 0 & A_3g_{n_3} & 0 & 0 & 0 \\
%  \dots & \dots & \dots & \dots & \dots & \dots \\
%  \dots & \dots & \dots & \dots & \dots & \dots \\
%  0 & 0 & 0 & 0 & 0 & A_kg_{n_k} \\
%\end{array}%
%\right),
%$$

%
%We are looking for a generalization where the local expression of
%$g_F$ is no longer block. As the base and each fiber are equipped
%with Riemannian (or semi-Riemannian) metrics, we can talk about
%local orthonormal sections of each fiber.

%Let $U_a \subset F_a$ be a trivializing
%neighborhood for $T^*F_a$, $V \subset B$ be a coordinate
%neighborhood on B and let \be W=U_1 \times U_2 \times ...\times
%U_k\times V.\ee Denote the local orthonormal sections of the
%cotangent bundle of each $F_a$ respectively by $\{\theta_a^i
%\}_{i=1}^{n_a}$ and the local coordinates on $B$ by
%$x^1,x^2,...,x^n$.
%
%A ``{\it warped-like product metric}" on $M$ is given by the local
%orthonormal frame of the cotangent bundle of $M$
%%
%\bea &&e_a^i=\sum_{b=1}^k \sum_{j=1}^{n_b} A_{aj}^{bi} \theta_b^j,
%\quad \ \ i=1,...,n_a, \ a=1,..,k \\ &&
% e_B^i=\sum_{j=1}^n a_{Bj}^i \ dx^j, \quad \ \ i=1,...,n \eea
% where \bea A_{aj}^{bi}= A_{aj}^{bi}(x_1,x_2,...,x_n),
%  a_{Bj}^{i}= a_{Bj}^{i}(x_1,x_2,...,x_n). \eea

\subsection{Special warped-like manifolds with 6-dimensional fiber space}

   {We will consider the case where the fiber $F$ is a
        $6$-manifold of the form $$ F=F_1\times F_2, $$ where
        $F_i$'s \ $(i=1,2)$ \ are $3$-manifolds each equipped
        with Riemannian metrics. Since all $3$-manifolds are
        paralellizable,  then the fiber $F$ is also
        paralellizable $6$-manifold. }

  {Let $\theta^i,\theta^{\hat{i}}$ \ $(i=1,2,3)$ \ be (global)
    orthonormal sections of the cotangent bundles of $F_1$ and
    $F_2$ respectively. Then the set of $3$-forms in the fiber
    $F$, i.e. $\Lambda^3(F)$, includes two closed $3$-forms.
    These are the volume forms of the $F_i$'s \ $(i=1,2)$
    given by $$ \label{volume1} vol_{F_1}=\theta^{123} \quad
    {\rm and} \quad vol_{F_2}=\theta^{\hat{1}\hat{2}\hat{3}}.
    $$}
  {We choose an {\it
almost complex structure } $J$ in the automorphism group of $TF$
 (i.e. $J\in End(TF)$ and $J^2=-I_d$) as follows,
\begin{eqnarray} \label{acs} J: & \ TF \quad \longrightarrow &TF \nonumber \\
&(\theta_i,\theta_{\hat{i}})\longmapsto
&(\theta_{\hat{i}},-\theta_{i}), \nonumber \end{eqnarray} where $\theta_i$ and
$\theta_{\hat{i}}$ are the dual of the $\theta^i$ and
$\theta^{\hat{i}}$ respectively. }

    An almost complex structure $J$ is called an {\it orthogonal
almost complex structure } if $J$ admits a Hermitian metric, $$
g(JX,JY)=g(X,Y), \ \ \forall X,Y \in \chi(F).$$  Note that every
almost complex structure admits a Hermitian metric on a
paracompact manifold
 and the Hermitian metric $h$ is given by $$
h(X,Y)=g(X,Y) +g(JX,JY).$$   Then the almost complex structure
admits a Hermitian metric.

 If we fix the Hermitian metric $g$ for the orthogonal almost
complex structure $J$, we obtain a non-degenerate  $2$-form
$\omega$ on $F$ and it is defined by $$ \omega(X,Y)=g(JX,Y),\ \
\forall X,Y \in \chi(F). $$ Note that the 2-form $\omega$ is also
$J$ invariant, i.e. $\omega(JX,JY)=\omega(X,Y)$.  Then we get
$\omega$ as $$ \label{kahler} \omega=\sum_{i=1}^3 \theta^i \wedge
J\theta^i=\theta^1\theta^{\hat{1}}+\theta^2\theta^{\hat{2}}+\theta^3\theta^{\hat{3}}.$$

   Using orthogonal almost complex structure $J$, it is
defined a complex volume form $\Psi$ on $F$ as $$ \label{complexv}
\Psi=\Psi^{+} +i
\Psi^-=(\theta^1+iJ\theta^1)(\theta^2+iJ\theta^2)(\theta^3+iJ\theta^3).
$$ Then we obtain \begin{eqnarray} \Psi^+&=& \theta^{123}
-\theta^{1\hat{2}\hat{3}} - \theta^{\hat{1}2\hat{3}}-\theta^{\hat{1}\hat{2}3}, \nonumber \\
\Psi^- &=& \theta^{\hat{1}23}
+\theta^{1\hat{2}3}+\theta^{12\hat{3}}-\theta^{\hat{1}\hat{2}\hat{3}}.
\nonumber \end{eqnarray}

  Hence the fiber space $F$ admits a special almost Hermitian $6$
dimensional manifold structure. Such a manifold is characterized
by its complex volume form $\Psi=\Psi^+ +i \Psi^-$ and its
K\"ahler form.

Let us consider a $(6+n)$-dimensional warped-like product manifold
$M=F_1 \times F_2 \times B$ where $F_1,F_2$ are 3-dimensional and
$B$ is an $n$-dimensional manifold. As we shall use this
structure, we
define it separately for ease of reference in the next studies. \\

\begin{definition} Let $M=F_1 \times F_2
 \times B$ be a $(6+n)$-dimensional topologically  product
 manifold   where $F_1$, $F_2$ are  $3$-manifolds and $B$ is
 an  $n$-manifold, each equipped with
 Riemannian metrics.
  Let $\theta^i,\theta^{\hat{i}}$ be orthonormal sections of the
cotangent bundles of $F_1$ and $F_2$ respectively and
$x_1,x_2,...,x_k$ be local coordinates on $B$. If  the metric on M
is defined by the following orthonormal frame
\bea \label{warped-like11} e^i&&=A(x_1,x_2,...,x_n)\theta^i + B(x_1,x_2,...,x_n) \theta^{\hat{i}}, \quad i=1,2,3 \nonumber \\
e^{\hat{i}}&&=\hat{A}(x_1,x_2,...,x_n)\theta^i + \hat{B}(x_1,x_2,...,x_n) \theta^{\hat{i}}, \quad i=1,2,3 \nonumber \\
e^{i+6}&&=a_{i1}(x_1,x_2,...,x_n)dx^1+...+a_{in}(x_1,x_2,...,x_k)dx^n,
\ \ \ i=1,2,...,n \eea where $A,B,\hat{A},\hat{B}$ and $a_{ij}$
$(i,j=1,2,...,n)$ are functions on base manifold $B$. Then we call
$(M,e^i)$ \ $(i=1,2,...,6+n)$ \ a {\it``$(3+3+n)$ warped-like
product"} manifold. \end{definition}

In the next section, we will consider a special case for the base
manifold $B$ is a two dimensional, that is, the fiber space $F$ is
a six dimensional manifold with $n=2$ dimensional base space.

\subsection{8-dimensional warped-like manifolds with 6-dimensional fiber space
%$F^6=F_1^3\times F_2^3$
}

The problem we are dealing with is modelled on an $8$-dimensional
warped-like product manifold $M=F_1 \times F_2 \times B$ where
$F_1,F_2$ are 3-dimensional and base manifold $B$ is a
2-dimensional manifold. Since we will use this structure often, we
define it   for ease of reference.

\begin{definition} \cite{BU}   Let $M=F_1 \times F_2
 \times B$ be an $8$-dimensional topologically  product
 manifold   where $F_1$, $F_2$ are  $3$-manifolds and $B$ is
 a  $2$-manifold, each equipped with
 Riemannian metrics.
  Let $\theta^i,\theta^{\hat{i}}$ be orthonormal sections of the
cotangent bundles of $F_1$ and $F_2$ respectively and $x$, $y$ be
local coordinates on $B$. If  the metric on M is defined by the
following orthonormal frame
\bea \label{warped-like} e^i&&=A(x,y)\theta^i + B(x,y) \theta^{\hat{i}}, \quad i=1,2,3 \nonumber \\
e^{\hat{i}}&&=\hat{A}(x,y)\theta^i + \hat{B}(x,y) \theta^{\hat{i}}, \quad i=1,2,3 \nonumber \\
e^{i+6}&&=a_{i1}(x,y)dx+a_{i2}(x,y)dy, \ \ \ i=1,2 \eea where
$A,B,\hat{A},\hat{B}$ and $a_{ij}$ $(i,j=1,2)$ are functions on
base manifold $B$. Then we call $(M,e^i)$ \ $(i=1,2,...,8)$ \ a
{\it``3+3+2 warped-like product"} manifold. \end{definition}

 The metric  given in the equation
(\ref{warped-like}) is written as \bea \label{ewl}
g&&=\sum_{i=1}^3 e^i \otimes e^i + e^{\hat{i}} \otimes e^{\hat{i}}
+ \sum_{i=1}^2 e^{i+6} \otimes e^{i+6} \nonumber \cr&&=
\underbrace{(a_{11}^2 + a_{21}^2) dx \otimes dx + (a_{12}^2 +
a_{22}^2 ) dy \otimes dy + ( a_{11}a_{12} + a_{22}a_{21}) dx
\otimes dy}_{\pi_B^* g_B} \nonumber \cr && +(A^2 + \hat{A}^2 )
\underbrace{\sum_{i=1}^3 \theta^i \otimes \theta^i}_{\pi_1^*
g_{F_1}} + (B^2 + \hat{B}^2 )\underbrace{\sum_{i=1}^3
\theta^{\hat{i}} \otimes \theta^{\hat{i}}}_{\pi_2^* g_{F_2}} +
2(AB+\hat{A} \hat{B}) \underbrace{\sum_{i=1}^3 \theta^i \otimes
\theta^{\hat{i}}}_{\pi_{12}^*\omega} \nonumber \\ \label{3219}&&
=\pi_B^* g_B + (f_1 \circ \pi_B) \pi_1^* g_{F_1} + (f_2 \circ
\pi_B) \pi_2^* g_{F_2} + (h_{12}\circ \pi_B) \pi_{12}^* \omega,
\eea where \be f_1=A^2 + \hat{A}^2, \quad f_2= B^2 + \hat{B}^2,
\quad h_{12}=2(AB+\hat{A} \hat{B}),\ee $\pi_B: F_1 \times F_2
\times
 \times B \longrightarrow B$, $\pi_{b}: F_1 \times
F_2  \times B \longrightarrow F_b$ and $\pi_{12}: F_1 \times F_2
 \times B \longrightarrow F_1 \times F_2$ are the
natural projections on $B$, $F_b$ and  $F_1 \times F_2$
respectively, and $\omega=\sum_{i=1}^3 \theta^i \theta^{\hat{i}}$
on $F_1\times F_2$. \\

\begin{remark} Note that $3+3+2$ warped-like  metric
is given as \cite{BU}   \be g_{3+3+2}=
 \left(%
\begin{array}{ccc|ccc|cc}
  A\theta^1 & 0 & 0 &  B\theta^{\hat{1}} &0 & 0 & 0  &0 \\
  0 & A\theta^2 & 0 & 0  & B\theta^{\hat{2}} & 0 & 0 & 0\\
  0 & 0 & A\theta^3 & 0  & 0 & B\theta^{\hat{3}} & 0& 0 \\ \hline
  \hat{A}\theta^1 & 0  & 0 &\hat{ B}\theta^{\hat{1}} & 0  & 0 &0& 0 \\
  0 & \hat{A}\theta^2 &0  & 0 & \hat{B}\theta^{\hat{2}} & 0 & 0& 0 \\
  0 & 0 & \hat{A}\theta^3 & 0  & 0 &\hat{B}\theta^{\hat{3}} & 0& 0 \\ \hline
  0 & 0 & 0 & 0 & 0 & 0 & dx & 0 \\ 0 & 0 & 0 & 0 & 0 & 0 & 0 & dy \\
\end{array}%
\right)\ee where $  \{ \theta^1,\theta^2,\theta^3,
\theta^{\hat{1}},\theta^{\hat{2}},\theta^{\hat{3}},dx,dy \} $ are
(global) frame on 8-dimensional manifold. \end{remark}

\subsection{Fibre-base decomposition on a (3+3+2) warped-like product manifold }

We present a special decomposition as an example of 8-dimensional
manifold (i.e. $8=3+3+2$). Corresponding to the decomposition of
the manifold as ``base" and ``fiber", the exterior algebra has the
following direct sum decomposition,
 \be
  \Lambda^p(M) =\bigoplus_{a+k=p} \Lambda^{(a,k)}(M),\ee
  where $a=1,\dots, 6$ and $k=1,2$, i.e. in our case the fiber is $6$-dimensional and the base is
$2$-dimensional, this fiber structure gives a decomposition of the
exterior algebra as follows.
\bea \Lambda^1(M)&=&\Lambda^{1,0}\oplus\Lambda^{0,1},\cr
\Lambda^2(M)&=&\Lambda^{2,0}\oplus\Lambda^{1,1}\oplus\Lambda^{0,2},\cr
\Lambda^3(M)&=&\Lambda^{3,0}\oplus\Lambda^{2,1}\oplus\Lambda^{1,2},\cr
\Lambda^4(M)&=&\Lambda^{4,0}\oplus\Lambda^{3,1}\oplus\Lambda^{2,2},\cr
\Lambda^5(M)&=&\Lambda^{5,0}\oplus\Lambda^{4,1}\oplus\Lambda^{3,2},\cr
\Lambda^6(M)&=&\Lambda^{6,0}\oplus\Lambda^{5,1}\oplus\Lambda^{4,2},\cr
\Lambda^7(M)&=&\Lambda^{6,1}\oplus\Lambda^{5,2},\cr
\Lambda^8(M)&=&\Lambda^{6,2}. \eea

Under the exterior derivative these summands are mapped as \be
d:\quad  \Lambda^{(a,k)}(M)\longrightarrow \Lambda^{(a+1,k)}\oplus
\Lambda^{(a,k+1)}.\ee We can refine this decomposition by
splitting the components for  each  fiber as
 \be
  \Lambda^p(M) =\bigoplus_{a+b+k=p}
\Lambda^{(a,b,k)}(M),  \ee where $a$ and $b$ range from $1$ to $3$
and $k=1,2$ as before, i.e., \beas
\Lambda^1(M)&=&\Lambda^{1,0,0}\oplus\Lambda^{0,1,0}\oplus\Lambda^{0,0,1},\cr
\Lambda^2(M)&=&\Lambda^{2,0,0}\oplus\Lambda^{1,1,0}\oplus\Lambda^{0,2,0}
\oplus\Lambda^{1,0,1}\oplus\Lambda^{0,1,1}\oplus\Lambda^{0,0,2},\cr
\Lambda^3(M)&=&\Lambda^{3,0,0}\oplus\Lambda^{2,1,0}\oplus\Lambda^{1,2,0}
\oplus\Lambda^{0,3,0}\oplus\Lambda^{2,0,1}\oplus\Lambda^{1,1,1}
\oplus\Lambda^{0,2,1}\oplus\Lambda^{1,0,2}\oplus\Lambda^{0,1,2},\cr
\Lambda^4(M)&=&\Lambda^{3,1,0}\oplus\Lambda^{2,2,0}\oplus\Lambda^{1,3,0}\oplus\Lambda^{3,0,1}
\oplus\Lambda^{2,1,1}\oplus\Lambda^{1,2,1}\oplus\Lambda^{0,3,1}\oplus\Lambda^{2,0,2}\oplus\Lambda^{1,1,2}
\cr && \oplus\Lambda^{0,2,2},\cr
\Lambda^5(M)&=&\Lambda^{3,2,0}\oplus\Lambda^{2,3,0}\oplus\Lambda^{3,1,1}
\oplus\Lambda^{2,2,1}\oplus\Lambda^{1,3,1}\oplus\Lambda^{3,0,2}\oplus\Lambda^{2,1,2}
\oplus\Lambda^{1,2,2}\oplus\Lambda^{0,3,2},\cr
\Lambda^6(M)&=&\Lambda^{3,3,0}\oplus\Lambda^{3,2,1}\oplus\Lambda^{2,3,1}\oplus\Lambda^{3,1,2}
\oplus\Lambda^{2,2,2}\oplus\Lambda^{1,3,2},\cr
\Lambda^7(M)&=&\Lambda^{3,3,1}\oplus\Lambda^{3,2,2}\oplus\Lambda^{2,3,2},\cr
\Lambda^8(M)&=&\Lambda^{3,3,2}.   \eeas

 The effect of the exterior derivative is
given by
 \be
 d: \quad  \Lambda^{(a,b,k)}(M)\longrightarrow
\Lambda^{(a+1,b,k)}\oplus \Lambda^{(a,b+1,k)}\oplus
\Lambda^{(a,b,k+1)}. \ee

%We shall first give a local  description of the
%generalization of the warped product structure.  We shall use the
%summation convention whenever appropriate.
%The base-fiber decomposition suggests that the local formulation of the
%line element is of the form
%$$ds^2=(g_F)_{ij}^{ab}\ dy_a^i\otimes dy_b^j +  (g_B)_{ij} dx^i\otimes dx^j,$$
%where $g_B$ depends on the base coordinates and $g_F$ depends both on the
%base and the fiber coordinates.
%For a  generalized warped product the local coordinate expression of $g_F$ is of the form
%$$g_F=\left(%
%\begin{array}{cccccc}
%  A_1g_{1} & 0 & 0 & 0 & 0 & 0 \\
%  0 & A_2g_{2} & 0 & 0 & 0 & 0 \\
%  0 & 0 & A_3g_{3} & 0 & 0 & 0 \\
%  \dots & \dots & \dots & \dots & \dots & \dots \\
%  \dots & \dots & \dots & \dots & \dots & \dots \\
%  0 & 0 & 0 & 0 & 0 & A_kg_{k} \\
%\end{array}%
%\right),
%$$
%where the $A_a$'s are functions of the base coordinates and each   $g_a$ the metric on $F_a$ hence
%is a symmetric matrix whose entries are functions of the coordinates of the fiber $F_a$.
The explicit example of $(3+3+2)$ warped-like product with
$Spin(7)$ manifold and $(3+3+1)$ warped-like product with $G_2$
manifold  will be presented in the next section (see details in
\cite{BU,U2,U}).

\section{Examples of warped-like product manifolds with $G_2$ and $Spin(7)$ holonomy}
 In this section, we present some special examples of warped-like product manifolds
    with $G_2$ and $Spin(7)$ holonomy group.

\subsection{(3+3+2) warped-like product manifold with $Spin(7)$ holonomy \cite{BU}}

We recall that the Yasui-Ootsuka solution \cite{YO} on \be
M=S^3\times S^3 \times {R}^2 \ee is given by the following
(global) orthonormal frame
 \bea \label{eeii222}
e^i&=&\frac{1}{2}{b}^{\frac{3}{4}}sech(y)\theta^i, \quad \quad
\quad \quad \quad \quad \quad \quad \
 i=1,2,3 \nonumber \\
e^{\hat{i}}&=&\frac{1}{2}ab^{-\frac{1}{4}}\left(1-tanh(y)\right)
\theta^i + ab^{-\frac{1}{4}} \theta^{\hat{i}}, \ \  i=1,2,3 \nonumber \\
e^7&=&a{b }^{\frac{3}{4}}dx, \nonumber \\  e^8 &=& {b
}^{\frac{3}{4}}sech(y) dy, \eea where the local
 sections of the
cotangent bundle of each $S^3$ respectively by $\theta^i,
\theta^{\hat{i}}$ and the functions $a(x)$, $b(x)$  satisfy the
differential equations \be \label{denk} \frac{da}{dx}
=\frac{1}{2}\left(\frac{a^3}{b}- a b \right),\quad
\frac{db}{dx}=-2a^2. \ee
 Thus the
metric is \bea \label{metrik} g&=& \underbrace{\left[ \sqrt{a^4b^3
} dx^2 + \sqrt{b^3 }sech^2(y) dy^2 \right]}_{\pi_B^* g_B}
\nonumber\\
 \label{2} &+&  \frac{1}{4}
\left[ \sqrt{b^3 }sech^2(y) + \sqrt{\frac{a^4}{b}} (1-tanh(y) )^2
\right]\underbrace{ \sum_{i=1}^3(\theta^i)^2}_{\pi_1^* g_{S^3}} +
\sqrt{\frac{a^4}{b}} \underbrace{\sum_{i=1}^3
(\theta^{\hat{i}})^2}_{\pi_2^* g_{S^3}}
\nonumber\\
&+& \label{3} \sqrt{\frac{a^4}{b}}
  \left[1-tanh(y)\right] \sum_{i=1}^3 \theta^i \theta^{\hat{i}},
 \eea
 where $\pi_B: S^3 \times S^3 \times {R}^2
\longrightarrow {R}^2$ and $\pi_i:  S^3 \times  S^3 \times {R}^2
\longrightarrow S^3$ are the natural projections on ${R}^2$ and
$S^3$ respectively.

Defining the functions \bea &&f_1= \frac{1}{4} \left[ \sqrt{b^3
}sech^2(y) + \sqrt{\frac{a^4}{b}} (1-tanh(y) )^2 \right], \quad
f_2= \sqrt{\frac{a^4}{b}}, \quad \nonumber \\ &&
h=\sqrt{\frac{a^4}{b}}
  [1-tanh(y)]
\eea and the 2-form $\omega$ \be \omega=\sum_{i=1}^3 \theta^i
\theta^{\hat{i}}, \ee we can write $g$ as \be \label{mulpwarp33}
g=\pi_B^* g_B +\sum_{i=1}^2(f_i\circ \pi_B)\pi_i^*
g_{F_i}+h\omega, \ee where $\pi_B: S^3 \times S^3 \times
{R}^2\longrightarrow {R}^2$ and $\pi_i:  S^3 \times S^3 \times R^2
\longrightarrow S^3$ are the natural projections on ${R}^2$ and
$S^3$ respectively. The matrix of $g$  with respect to the
(global) frame \be \{
\theta^1,\theta^2,\theta^3,\theta^{\hat{1}},\theta^{\hat{2}},\theta^{\hat{3}},e^7,e^8
\} \ee is \be g=
\left(%
\begin{array}{ccc}
  f_1I_3 & \frac{h}{2}I_3       & 0   \\
    \frac{h}{2}I_3    & f_2I_3  & 0   \\
    0    & 0       & I_2 \\
\end{array}%
\right) \ee where $I_3$ and $I_2$ are identity matrices of size
$3$ and $2$, and zeroes denote zero matrices of appropriate sizes.
Note that if $h$ were zero, the metric given in the equation
(\ref{mulpwarp33}) would be a multiply warped product with a block
diagonal matrix with respect to an appropriate frame \cite{FS}.

 We present a parallel
$Spin(7)$ manifold as a special warped-like product structure on
\be M=S^3\times S^3 \times R^2 \ee \cite{YO,BU}. If we choose the
(global) orthonormal frame   as in the equation (\ref{baz}),
 then the coefficient (in the sense of warping) functions are \cite{BU},
 \bea &&
A=\frac{1}{2}{b}^{\frac{3}{4}}sech(y), \quad B=0, \quad
\hat{A}=\frac{1}{2}ab^{-\frac{1}{4}}(1-tanh(y)),  \quad
\hat{B}=ab^{-\frac{1}{4}}, \nonumber \\&& a_{11}=a{b
}^{\frac{3}{4}},\quad a_{12}=0,\quad a_{21}=0,\quad a_{22}={b
}^{\frac{3}{4}}sech(y).\eea  Hence the metric on $M$ is given by
\bea \label{metrik} g&&= \left[ \sqrt{a^4b^3 } dx^2 + \sqrt{b^3
}sech^2(y) dy^2 \right]
 \label{2}  +   \frac{1}{4}
\left[ \sqrt{b^3 }sech^2(y) + \sqrt{\frac{a^4}{b}} (1-tanh(y) )^2
\right]  \sum_{i=1}^3(\theta^i)^2  \nonumber\\
&&+ \sqrt{\frac{a^4}{b}}
 \sum_{i=1}^3 (\theta^{\hat{i}})^2
+ \label{3} \ \ \ \sqrt{\frac{a^4}{b}}
  \left[1-tanh(y)\right] \sum_{i=1}^3 \theta^i \theta^{\hat{i}}, %\nonumber
 \eea where the local
 sections of the
cotangent bundle of each $S^3$ respectively by $\theta^i,
\theta^{\hat{i}}$ and the functions $a(x)$, $b(x)$  satisfy the
differential equations \be \label{denk} \frac{da}{dx}
=\frac{1}{2}\left(\frac{a^3}{b}- a b \right),\quad
\frac{db}{dx}=-2a^2. \ee $\hfill \square $ \\

The following theorem is an important result for Yasui-Ootsuka
solution and (3+3+2) warped-like product manifold with $Spin(7)$
holonomy. See details in \cite{BU,BU2}

\begin{theorem} \cite{BU}     Let M be diffeomorphic
to $F \times B$, where the base B is a two-dimensional Riemannian
manifold diffeomorphic to $R^2$ , the fiber F is a 6-manifold of
the form $F = F_1 \times F_2$, and the $F_i$ (i = 1, 2) are
complete, connected and simply connected 3-manifolds. Let the
metric on M be a (3 + 3 + 2) warped-like product. Then the fibers
$F_i$ are isometric to $S^3$ with constant curvature $k > 0$ for
the Spin(7) structure determined by the Bonan form (5). Also there
exists a unique metric in the class of (3 + 3 + 2) warped-like
product metrics admitting the same Spin(7) structure, and the
metric is written as given by Yasui-Ootsuka (\ref{metrik}) up to a
gauge transformation.  \end{theorem}

\subsection{(3+3+1) warped-like product manifold with $G_2$ holonomy \cite{U}}

Inspired by the previous work \cite{BU}, moving down one
dimension, we study $(3+3+1)$ warped-like product manifolds
\cite{U}. In \cite{BU}, we studied special type of warped-like
product and its properties with $Spin(7)$ holonomy in
$8$-dimensional manifolds. Here the problem we are dealing with is
modeled on an $7$-dimensional warped-like product manifold $M=F_1
\times F_2 \times B$ where $F_1,F_2$ are 3-dimensional and $B$ is
a one dimensional manifold.

 Using by the  paper \cite{BU}, we
will define $(3+3+1)$ warped-like product manifolds  in this
section  as an example of $G_2$ manifold \cite{U}.

\begin{definition} Let $M=F_1 \times F_2
 \times B$ be an $7$-dimensional topologically  product
 manifold   where $F_1$, $F_2$ are  $3$-manifolds and $B$ is
 a  one dimensional manifold, each equipped with
 Riemannian metrics.
  Let $\theta^i,\theta^{\hat{i}}$ be orthonormal sections of the
cotangent bundles of $F_1$ and $F_2$ respectively and $x$ be local
coordinate on $B$. If  the metric on M is defined by the following
orthonormal frame \bea  e^i&=&A(x)\theta^i + B(x)
\theta^{\hat{i}}, \nonumber\\ e^{\hat{i}}&=&\hat{A}(x)\theta^i +
\hat{B}(x) \theta^{\hat{i}}, \label{baz} \\  e^{7}&=&a(x)dx, \quad
i=1,2,3 \ \ \nonumber \eea then we call $(M,e^i)$ \ $i=1,2,...,7$
a {\it``(3+3+1) warped-like product"} manifold. \end{definition}

\subsubsection{An example of (3+3+1) warped-like manifold with $G_2$ holonomy} Recall
that the Konishi-Naka solution \cite{Naka} on \be M=SU(2)\times
SU(2) \times {R}  \ee is given by the following (global)
orthonormal frame
 \bea \label{eeii222}
e^i&=&A(x)\theta^i, \quad \quad \quad \quad \quad \quad \ \
 i=1,2,3 \nonumber \\
e^{\hat{i}}&=&\widehat{A} \left(
   \theta^{\hat{i}}-\frac{1}{2}\theta^i \right), \   \quad \quad \quad  i=1,2,3 \nonumber \\
e^7&=&dx,  \eea where the local
 sections of the
cotangent bundle of each $SU(2)$ respectively by $\theta^i,
\theta^{\hat{i}}$ and the functions $A(x)$, $\widehat{A}$  satisfy
the differential equations \be \label{denkkk} \frac{d
 {A}}{dx} =  \frac{\widehat{A}}{2A},\quad \frac{d\widehat{A}}{dx}=1-\frac{ \widehat{A}^2}{4A^2}. \ee

 Thus the
metric is \bea \label{KNmetrik} g= A(x)^2\sum_{i=1}^3(\theta^i)^2
+ \widehat{A}(x)^2\left(\sum_{i=1}^3 \left[\theta^{\hat{i}}
-\frac{1}{2}  \theta^i \right]  \right)^2 +dx^2. \eea

 We can
take $ e^7=dx,  $ as in \cite{Naka}. It is seen that this metric
is an example of (3+3+1) warped-like product manifold with $G_2$
holonomy. See details in \cite{U,Naka}. $\hfill \square$ \\

The following result is important for Konishi-Naka solution and
(3+3+1) warped-like product manifold with $G_2$ holonomy. \\

\begin{theorem}
 \cite{U}   Let $M$ be diffeomorphic
to $F \times B$, where the base $B$ is a one dimensional
Riemannian manifold diffeomorphic to $R$, the fibre $F$  is a
$6$-manifold of the form $F=F_1\times F_2$, and $F_i$ $(i=1,2)$
are complete, connected and simply connected $3$-manifolds. Let
the metric on M be a   (3+3+1) warped-like product metric. If $M$
is the manifold with the $G_2$  holonomy or with the weak $G_2$
holonomy determined by the fundamental 3-form (1), then the fibers
$F_i$'s are isometric to $S^3$ with constant curvature $k>0$.
% Let the metric on $M$ be a special warped-like
% product. Then the fibers $F_i$'s are isometric to $S^3$ with constant curvature $k>0$
%for the $G_2$ holonomy and weak $G_2$ structure determined by the
%fundamental 3-form (\ref{form}).
  Also there
exists a unique metric in the class of special warped-like product
metrics admitting the   $G_2$ holonomy, and the metric is written
as given  by Konishi-Naka (\ref{KNmetrik}) up to gauge
transformation.
\end{theorem} %\hfill $\square$

In the next sections, we present some other explicit examples of
special warped-like metric constructions.

\subsection{Special warped-like product manifold with $G_2$ holonomy \cite{Bra}}

The  general metric ansatz compatible with  $G_2$ holonomy is
given \cite{Bra} by \be ds^2= \sum_{a=1}^7 e^a \otimes e^a \ee
with the following  \bea  &&e^1 = A(r)(\sigma^1 - \Sigma^1 ),  \
e^2 = A(r)(\sigma^2 - \Sigma^2 ),  \nonumber \\ && e^3 =
D(r)(\sigma^3 - \Sigma^3 ),  \ e^4 = B(r)(\sigma^1 + \Sigma^1 ),
\nonumber \\ && e^5 =B(r)(\sigma^2 + \Sigma^2 ),  \ e^6 =
C(r)(\sigma^3 + \Sigma^3 ),  \nonumber \\ && e^7 =
\frac{dr}{C(r)}, \eea where $\sigma^i$ and $\Sigma^i$ \
($i=1,2,3$) are two sets of left invariant $SU(2)$ one-forms (see
details in \cite{Bra}) and their metric have made a particular
choice for the radial coordinate. Their ansatz depends on four
functions. Note that this special warped-like  metric is given as
\cite{Bra} be \be g_{special}=
 \left(%
\begin{array}{ccc|ccc|cc}
  A\sigma^1 & 0 & 0 &  -A\Sigma^{ {1}} &0 & 0 & 0   \\
  0 & A\sigma^2 & 0 & 0  & -A\Sigma^{ {2}} & 0 & 0 \\
  0 & 0 & D\sigma^3 & 0  & 0 & -D\Sigma^{ {3}} & 0\\ \hline
  B\sigma^1 & 0  & 0 &B\Sigma^{ {1}} & 0  & 0 &0 \\
  0 & B\sigma^2 &0  & 0 & B\Sigma^{ {2}} & 0 & 0 \\
  0 & 0 & C\sigma^3 & 0  & 0 &C\Sigma^{ {3}} & 0 \\ \hline 0 & 0 & 0 & 0 & 0 & 0  & dr \\
\end{array}%
\right)\ee where $  \{ \theta^1,\theta^2,\theta^3,
\theta^{\hat{1}},\theta^{\hat{2}},\theta^{\hat{3}},dr \} $ are
(global) frame on 7-dimensional manifold. To prove that the above
metric has $G_2$ holonomy, we have to impose that the associative
three-form $\varphi$ is closed and co-closed (i.e.
$d\varphi=d*\varphi=0$). These conditions imposed on the
associative three-form $\varphi$ constructed from metric structure
leads to the following system of first order differential
equations \bea \frac{dA}{dr} &=&
\frac{1}{4}\left(\frac{B^2-A^2+D^2}{BCD} +\frac{1}{A} \right)
 \nonumber \\  \frac{dB}{dr} &=&
\frac{1}{4}\left(\frac{A^2-B^2 +D^2}{ACD} -\frac{1}{B} \right)
 \nonumber\\   \frac{dC}{dr} &=&
\frac{1}{4}\left(\frac{C}{B^2} -\frac{C}{A^2} \right)
 \nonumber\\    \frac{dD}{dr} &=&
\frac{1}{2}\left(\frac{A^2+B^2 -D^2}{ABC}   \right)
   \eea
Special solutions is presented in their paper \cite{Bra}. This
metric structure present an example of special warped-like product
manifold with $G_2$ holonomy.

%
%\subsection{Special warped-like product manifold with $G_2$ holonomy \cite{Cve}}
%
%From the spin connection for the metric, we find that in the
%obvious choice of orthonormal basis
%
%\beas  && e^0 = dt, e^i = b  \sigma^i \\ && e^{\tilde{i}} = a (
%\Sigma^i+(g-\frac{1}{2} ) \sigma^i)\eeas
%
%
%

\subsection{Another special warped-like product manifold with $G_2$ holonomy \cite{Cve0}}
In this section we present  $G_2$ metrics for the six-function
$a_i,b_i \ (i=1,2,3)$ on $S^3 \times S^3$ manifold. A rather
general ansatz involving more functions was considered in
\cite{Bra,Cve}, and first-order equations for $G_2$ holonomy were
derived. The metric for the six-function $G_2$ space is given by
\be ds^2= dt^2 + a_i^2(\sigma^i- \Sigma^i) + b_i^2 (\sigma^i +
\Sigma^i)\ee where $\sigma^i$ and $\Sigma^i$ are left-invariant
1-forms for two $SU(2)$ group manifolds. It is seen that this
special warped-like metric is given as   be \be g_{special}=
 \left(%
\begin{array}{ccc|ccc|cc}
  a_1\sigma^1 & 0 & 0 &  -a_1\Sigma^{ {1}} &0 & 0 & 0   \\
  0 &a_2\sigma^2 & 0 & 0  & -a_2\Sigma^{ {2}} & 0 & 0 \\
  0 & 0 & a_3\sigma^3 & 0  & 0 & -a_3\Sigma^{ {3}} & 0\\ \hline
  b_1\sigma^1 & 0  & 0 &-b_1\Sigma^{ {1}} & 0  & 0 &0 \\
  0 & b_2\sigma^2 &0  & 0 & -b_2\Sigma^{ {2}} & 0 & 0 \\
  0 & 0 & b_3\sigma^3 & 0  & 0 &-b_3\Sigma^{ {3}} & 0 \\ \hline 0 & 0 & 0 & 0 & 0 & 0  & dt \\
\end{array}%
\right)\ee where $  \{ \theta^1,\theta^2,\theta^3,
\theta^{\hat{1}},\theta^{\hat{2}},\theta^{\hat{3}},dt \} $ are
(global) basis on 7-dimensional manifold. It was found that for
$G_2$ holonomy, $a_i$ and $b_i$ must satisfy the first-order
equations

\beas \frac{da_1}{dt} &=& \frac{a_1^2}{4a_3b_2}+
\frac{a_1^2}{4a_2b_3}-\frac{a_2}{4
b_3}-\frac{a_3}{4b_2}-\frac{b_2}{4a_3}-\frac{b_3}{4a_2} \\
\frac{da_2}{dt} &=& \frac{a_1^2}{4a_3b_2}+
\frac{a_1^2}{4a_2b_3}-\frac{a_2}{4
b_3}-\frac{a_3}{4b_2}-\frac{b_2}{4a_3}-\frac{b_3}{4a_2} \\
\frac{da_3}{dt} &=& \frac{a_1^2}{4a_3b_2}+
\frac{a_1^2}{4a_2b_3}-\frac{a_2}{4
b_3}-\frac{a_3}{4b_2}-\frac{b_2}{4a_3}-\frac{b_3}{4a_2}
\\ \frac{db_1}{dt} &=& \frac{a_1^2}{4a_3b_2}+
\frac{a_1^2}{4a_2b_3}-\frac{a_2}{4
b_3}-\frac{a_3}{4b_2}-\frac{b_2}{4a_3}-\frac{b_3}{4a_2}
\\ \frac{db_2}{dt} &=& \frac{a_1^2}{4a_3b_2}+
\frac{a_1^2}{4a_2b_3}-\frac{a_2}{4
b_3}-\frac{a_3}{4b_2}-\frac{b_2}{4a_3}-\frac{b_3}{4a_2} \\
\frac{db_3}{dt} &=& \frac{a_1^2}{4a_3b_2}+
\frac{a_1^2}{4a_2b_3}-\frac{a_2}{4
b_3}-\frac{a_3}{4b_2}-\frac{b_2}{4a_3}-\frac{b_3}{4a_2}. \eeas

Under some special conditions and numerical analysis on this
system given above, a special solution of these equations is
presented in \cite{Cve0}. \hfill $\square$ \\

\begin{remark} One can consider the relation cone metric structure
with special warped-like product manifold with $G_2$ (or
$Spin(7)$) holonomy. Consider the cone $M=F \times R$ with the
metric \be \label{321} ds^2= dt^2 + A_i^2(\sigma^i+ \Sigma^i) +
B_i^2 (\sigma^i - \Sigma^i)\ee where $\sigma^i$ and $\Sigma^i$ are
the standard coframe of 1-forms, $ A_i(t)$ and $B_i(t)$ are
positive functions defining the deformation of the standard cone
metric on M (see details in \cite{Baza}).
\end{remark}

 Complete Riemannian
metrics with holonomy group $G_2$ on manifolds obtained by
deformation of cones over $S^3 \times S^3$ are constructed in
\cite{Baza}. Their   idea of the paper is to consider the standard
conic metric over a Riemannian manifold with special geometry.
After then any deformation of this metric depends on a number of
functional parameters using which the $G_2 \ (or \ Spin(7))$
structure can be explicitly specified. In \cite{Baza},  they
propose (regarding to \cite{Bra,Cve0}) to consider the fiber space
$F=S^3 \times S^3$. Then the cone metric can be written as  in
equation \eqref{321}. A system of differential equations
guaranteeing that the metric $ds^2$ has a holonomy group contained
in $G_2$ was written out in \cite{Bra}. A particular solution of
this system corresponding to the metric with holonomy group $G_2$
on $S^3\times  R^4$ was also obtained in \cite{Bra}. Note that
more general metrics on cones over $S^3 \times S^3$ were also
studied in \cite{Bra,Cve0} (see also refs. in \cite{Baza}). In the
paper \cite{Baza}, they continue the study of this class of
metrics by setting the particular case $A_2=A_3$ and $B_2=B_3$ and
considering a boundary condition different from the one in
\cite{Bra}, which leads to spaces with
other topological structures. %Namely, we require that only the
%function B1 vanish at the vertex of the cone.

\begin{remark}
It is also noted that the metrics of Bryant and Salamon \cite{RS}
are asymptotically conical. In other words, outside of
sufficiently large compact sets, the metric  is arbitrarily close
to a conical metric $g_t=dt^2 + t^2 g$ defined on $M=S^3 \times
S^3\times R^+$. The closeness is measured with respect to the
conical metric \cite{Clarke}. A consequence of $(M, g_t )$ having
$G_2$ holonomy is that $(S^3 \times S^3, g )$ must be nearly
Kahler (see details in \cite{RS,Clarke}). In each given case they
were able to solve a system of ordinary differential equations so
as to obtain a metric $g_t$ of exceptional holonomy.

\end{remark}

\section{Conclusions}

In this paper we define warped-like product metrics as a
generalization of multiply-warped products and present examples of
special type of these metrics for $G_2$ and $Spin(7)$ holonomies
(i.e. exceptional holonomy cases \cite{Bra,Cve0,Jb,Salamonb}). We
investigate a generalization of the warped product as a special
warped-like form where the local coordinate expression of the
metric is no longer block diagonal, but still with a certain
simple structure. Since we want to have a fiber structure, we
shall allow dependencies on the coordinates of the base
everywhere. To obtain more explicit examples, different types of
fiber-base decompositions will be worked in the next studies. The
study of special holonomy manifolds (i.e. Berger's list \cite{Jb})
and some classes of $G_2$, $Spin(7)$ (see Tables 1-2) manifolds
with warped-like metric structure are  open problems for future
research topics. Finally we believe that our approach of the
warped-like product metrics will be an important notion  on the
manifolds with special and exceptional holonomies for future
studies, and also other types of geometries which use warped and
multiply-warped product.

\end{document}